# Student-Centric Graduate Training in Mathematics:

# A Commentary

Yara Skaf and Reinhard Laubenbacher
University of Florida

Yara.skaf@ufl.edu
Reinhard.laubenbacher@medicine.ufl.edu

**Abstract**

Career opportunities for PhDs in the mathematical sciences have never been better. Traditional faculty positions in mathematics departments in colleges and universities range from all teaching to combined teaching/research responsibilities. Beyond those, a wide array of careers has now opened up to freshly minted graduates, in academics, industry, business, and government. It is well-understood that these all require somewhat different preparations for Ph.D.s to be competitive. This commentary compares and contrasts mathematics graduate programs with Ph.D. programs in the life and biomedical sciences, which are structured in a way that allows considerable customization around students' career goals. While these programs may not be appropriate templates for the mathematical sciences, they have some features that might be informative. This commentary is intended to add perspective to the ongoing discussion around PhD training in the mathematical sciences. It also provides some concrete proposals for changes.

**Introduction**

In the 2018 National Academies report "Graduate STEM Education for the 21st Century," one of the many recommendations for reform of U.S. graduate education is to make programs more student-centric, taking into account their career goals and interests, coupled with extensive advising. For the mathematical sciences, this is timely advice, as the professional landscape for today's mathematical scientists is rapidly evolving. Compared to 7.4% across all occupations, employment in mathematics and statistics is projected to increase by 27.9% from 2016 to 2026, with most opportunities related to novel applications of mathematics and computational sciences to emerging fields such as digital cryptography, financial modeling, and healthcare analytics. This coincides with what seems to be a widespread sentiment among students and faculty alike that employment prospects for mathematics PhDs in academia are declining-- many graduate students appear skeptical about the likelihood of obtaining a tenured faculty position. As several studies have shown, there are many more PhDs in STEM fields than there are positions in academia, and this gap is only getting wider. Further, surveys across STEM disciplines indicate that the prospect of continual stress and uncertainty make the academic track distinctly unappealing for many students.[11] Taken together, these trends make for a compelling argument for PhD program reform. In this commentary we want to draw attention to features of some other program templates that can provide a different perspective. We chose programs in the life and biomedical sciences for comparison, in particular MD/PhD programs that prepare students to be both physicians and scientists, as these are the non-mathematics programs with which we are most familiar.

The traditional model of doctoral education in mathematics is effective at preparing trainees for careers in academic departments with a dual research/teaching mission via the standard trajectory of graduate school followed by postdoctoral training and eventually teaching/research as tenure-track faculty. But, as mentioned, many PhD students will ultimately pursue career paths distinct from this trajectory. So programs that acknowledge this reality and adapt their training accordingly will serve their students better and allow them to fully take advantage of the myriad exciting career opportunities inside and outside of academia. This discussion around effective training of doctoral students for modern careers is of course by no means new-- a robust body of literature seeks to analyze and propose reform in graduate education.[1-9] Here, we hope to contribute a unique perspective on this issue, stemming from our own experiences with various types of training programs at several different institutions. We first examine and compare the structures of current doctoral programs across some STEM fields as well as the career trajectories of their graduates. Based on this assessment, we identify a number of concerns related to graduate training in mathematics. Finally, we propose a number of modifications to potentially mitigate these concerns. Our aim is to highlight areas where standard practices could continue to adapt to better prepare students for their eventual careers.

**Graduate Program Structure**

A survey of the literature on graduate education in STEM fields resulted in a number of observations that affirm our own personal experiences. First, most programs take what is referred to as a "one-size-fits-all" approach, meaning the program elements and requirements are essentially the same for every student. Of the math graduate programs ranked in the top 50 by US News & World Report,[23] very few make any explicit mention of options to customize degree requirements to better suit individual student needs. Even these rare opportunities for customization are limited to things like fewer departmental exams in pure math for applied math students and occasional options to pursue industry fellowships instead of teaching assistantships. Second, the interests and career trajectories of graduate trainees extend far beyond tenured positions in academic research. Finally, despite similar trends toward diversification of professional opportunities and eventual employment choices for graduates (see Figure 1), mathematics programs remain structurally distinct from their counterparts in the life and biomedical sciences, in particular with regard to their strong focus on preparing students for careers in academia. Virtually all of the top 50 programs mention preparing students for careers in research or academia, but far fewer mention careers in industry or government, and fewer still explicitly offer any program elements designed to facilitate these alternative career paths.

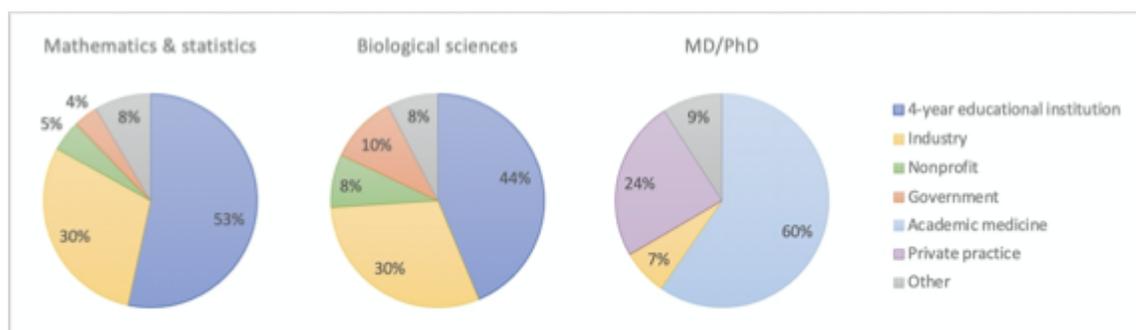

**Figure 1.** Employment sector of employed PhD and MD/PhD graduates in mathematics and biological sciences. Note that the proportion of doctoral researchers employed by academic institutions (or in academic medicine for MD/PhDs) is relatively similar in all three populations.

The vast majority of graduate programs divide their focus between three major areas: (1) coursework, which provides a rigorous theoretical foundation for work in the field of study; (2) research activities, which provide practical training in independent investigation, scientific communication, and grantsmanship; and (3) teaching responsibilities, which provide financial support throughout the program, in addition to experience with mentorship and pedagogy. Most programs also include (4) departmental exams, though these are generally taken early in training so that this component only appears in the initial years. Despite having identical components, the relative emphasis on each element differs significantly across domains. (Table 1) For instance, mathematics programs tend to include substantial coursework and teaching while PhDs in the biomedical sciences generally focus primarily on research.

**Mathematics.** In addition to dissertation research, standard PhD requirements consist of preliminary exams in several core subjects, oral and written qualifying exams in a chosen focus area, and substantial coursework throughout the program. Since stipends are mostly funded through teaching, graduate students must tutor, grade for, and assist with multiple undergraduate classes every semester or teach an entire class themselves.

**Biomedical sciences.** By contrast, most life science PhDs require minimal coursework and no preliminary exams. Any mandatory classes are usually closer to seminars or journal clubs in which students read and discuss research papers, but these contain little to no new theory beyond what was learned in an undergraduate curriculum. Compared to 2-3 intensive graduate mathematics courses with exams and assignments that are expected to take many hours every week, participating in 1-2 seminar-style classes each semester requires significantly less time. Additionally, many graduate students assist with labs instead of lectures, where their main role is supervision of students who independently perform experiments. Most undergraduate biology lecture courses give only multiple choice exams, so teaching assistants spend minimal time on grading. Though qualifying examinations are required, the written portion generally takes the form of a drafted research proposal rather than the lengthy theory-based exams given in mathematics

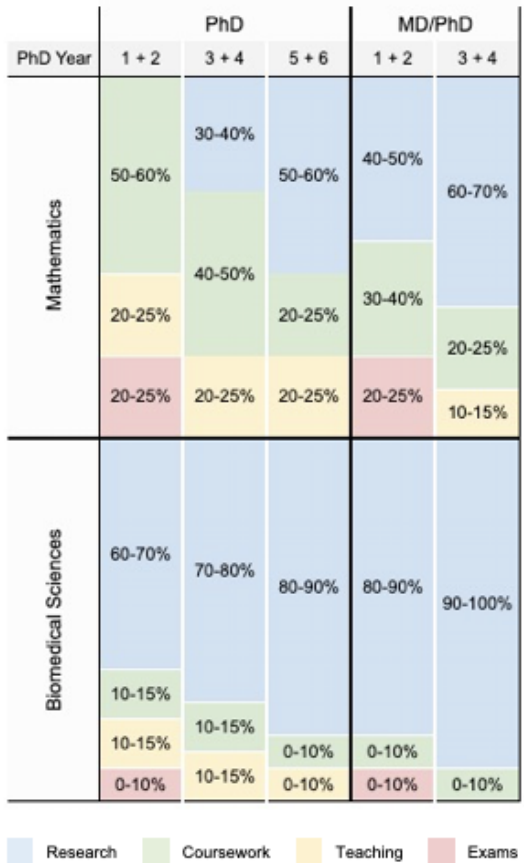

**Table 1.** Composition of typical graduate programs.

| PhD Year | PhD 1 + 2 | PhD 3 + 4 | PhD 5 + 6 | MD/PhD 1 + 2 | MD/PhD 3 + 4 |
|---|---|---|---|---|---|
| **Mathematics** | 50-60% / 20-25% / 20-25% | 30-40% / 40-50% / 20-25% / 20-25% | 50-60% / 20-25% / 20-25% | 40-50% / 30-40% / 20-25% | 60-70% / 20-25% / 10-15% |
| **Biomedical Sciences** | 60-70% / 10-15% / 10-15% / 0-10% | 70-80% / 10-15% / 10-15% | 80-90% / 0-10% | 80-90% / 0-10% / 0-10% | 90-100% / 0-10% |

Legend: Research | Coursework | Teaching | Exams

programs. These savings from decreased time on coursework, exams, and teaching are instead invested in research-- students immediately begin rotating through different faculty labs to select an advisor for their PhD project, which should be well underway by the end of the first year. Students are also expected to assist with grant proposals, prepare manuscripts, and present original research with a much higher degree of regularity than would be typical for mathematics students.

**MD/PhD programs.** The primary focus of MD/PhD programs is training students to conduct translational research generating or applying basic science results to improve medical care. As such, the PhD portion of these programs tends to emphasize research activities almost exclusively. A typical MD/PhD program includes four years of biomedical science PhD training, interspersed with four years of medical school. To enable this accelerated completion of a PhD in four years instead of six, most programs have very few requirements. Medical school classes are allowed to fulfill the vast majority of Graduate School course requirements, and any mandatory courses are generally of the seminar form described above. Similarly, there are no preliminary exams, and qualifying exams follow the biomedical sciences format.

MD/PhD trainees are funded directly through their program and/or federal grants. This eliminates the need for students to generate their own stipends by teaching, and thus allows even more time for research than PhD-only programs. Continuous involvement in research is a firm expectation from day one onward, with many institutions encouraging students to begin graduate work as early as the summer before the official program start date. This expectation includes substantive contributions to all aspects of the investigative process from planning and funding acquisition to execution, presentation, and publication.

**Variability across programs.**

This comparison serves to illustrate alternative approaches to graduate education. Many factors undoubtedly contribute to this structural contrast, in particular related to variation in career trajectories, measures of academic success, and intrinsic differences between the underlying domains. We discuss these differences here in some detail.

*Undergraduate education.* Some interdisciplinary variation in program structure is attributable to differences in the preparedness of entering PhD students to independently conduct research, which is rooted in their academic backgrounds before graduate school. Most first year graduate students in biology are immediately capable of driving and executing significant research projects on their own. One reason for this is that the theoretical foundation for biology is constructed long before graduate school-- high school courses convey the majority of core concepts, with undergraduate courses reinforcing these ideas while adding some additional level of detail. Toward the end of their undergraduate program, students have accrued enough basic knowledge for classes to shift beyond theoretical foundations toward practical discussions of journal articles. Early PhD students have thus already developed the ability to critically assess research at the cutting edge of their field. A second reason is that biology students can gain experience and develop practical lab skills before they necessarily have a firm grasp of the underlying theory— an inexperienced researcher can learn how to follow a protocol for genomic sequencing library preparation years before they understand the biology underlying those tools.

While beginning graduate students in biology typically have substantial theoretical background, laboratory skills, and research experience that allow them to function largely independently in the lab, mathematics students usually require several years of training before they have the necessary skills for research. Unlike the life sciences, high school coursework rarely progresses beyond basic calculus, meaning most if not all theoretical mathematics is learned during their undergraduate programs. At the beginning of graduate school, the majority of students have foundational knowledge in the core areas of mathematics, but have yet to encounter the level of advanced theory that is needed to tackle a research problem. Since there is no real way to read or conduct research in mathematics without a solid understanding of the principles underpinning the results, it takes students longer to reach a point where they can make meaningful contributions to the field. Graduate programs in mathematics must therefore spend more time on coursework in core knowledge areas.

*Professional expectations.* Another contributor is the difference in professional expectations between a mathematics PhD graduate and a life sciences PhD or MD/PhD graduate. Graduate programs in any field generally have as a primary goal the training of students for positions as research faculty in academia. In the biomedical sciences, this means developing competency with grant applications, publishing papers, presenting work at conferences, and managing less experienced scientists in addition to designing and executing research initiatives. As a result, each of these activities is a substantial component of any typical biology PhD curriculum. By contrast, mathematics PhD students are expected to spend most of their first two years taking courses, teaching, and preparing for departmental examinations. Graduates are expected to complete enough research to publish their dissertation in a professional journal, then pursue postdoctoral or other entry-level positions where they can gain experience with independent research. There is somewhat less emphasis in mathematics on acquiring grant funding and generating a constant stream of publishable work. Instead, faculty are expected to teach and mentor students as they pursue longer-term research. As such, programs focus on helping students develop solid theoretical foundations and teaching experience.

This difference is even more pronounced for MD/PhDs. Essentially all MD/PhD students have been involved with research throughout their undergraduate education—substantial research experience is a requirement for admission to any major program. Students are therefore capable of making progress on graduate research immediately upon starting the program. During the PhD years, most publish their own research articles and carry out projects independently. Since some financial burden on the institution can be alleviated if students obtain their own funding, they are typically encouraged or required to independently apply for grants and fellowships to offset the cost of tuition and stipends whenever possible. Upon completion of their PhD requirements, students return to medical school, go on to a clinical residency program, and ultimately pursue clinical or translational research faculty positions in academic medicine. Clinical residencies allow minimal time for research, so at the time of graduation, MD/PhDs are already expected to be virtually fully independent physician scientists who can contribute funding and novel research to their future institution.

Putting aside the reasons for these differences between programs, we see (Table 1) that, compared to mathematics PhD programs, other programs spend proportionally more time on student research and, as a result, the programs can be easily customized to students' interests and career aspirations. One interesting difference to observe is that graduate students in biomedical programs come better prepared to carry out research projects early on in their career. In comparison, undergraduate mathematics programs typically have limited research opportunities for students beyond that provided by NSF-funded Research Experiences for Undergraduates programs. Another noteworthy difference is that graduate courses in mathematics are typically more heavily lecture- and exam-based. Qualifying examinations consist of extensive technical skill testing in mathematics versus research proposal writing in biomedical programs.

### Discussion

**Concerns with the current mathematics model.** There are a number of issues associated with the current "one-size-fits-all" structure of mathematics PhD programs. First and foremost, this structure does not optimally prepare many students for the myriad of diverse career paths in mathematics, especially in applied mathematics. Despite only about half of doctoral mathematicians being employed in academia (far fewer in tenured positions), graduate education remains focused on this professional trajectory. Such a program composition makes it challenging to find time for research, especially during the initial years where coursework and teaching are the most demanding. As a result, students graduate with comparatively less research experience than their colleagues in other fields. On the opposite end of the spectrum, the time spent on dissertation work in later years provides minimal benefit to students who plan to pursue careers as lecturers or other positions that involve little to no research.

**Examples of alternatives.** We mention two concrete examples of PhD programs that contain elements worth studying. One extreme example of a highly flexible PhD program is at the Okinawa Institute of Science and Technology (OIST). The PhD program at this graduate university is individualized to each student and designed to provide significant training and opportunity for cross-disciplinary research in addition to foundational knowledge about a primary field of interest.[17] To encourage cross-disciplinary collaboration, OIST has no formal departments-- research faculty are partitioned into "units" according to their topic of study.[18] As a testament to the research success of this model, OIST was ranked first in Japan and ninth worldwide for the

proportion of research papers published in leading scientific journals by the Nature Index[19] in 2018.[20]

A somewhat less extreme example is the interdepartmental PhD program in computational biology and bioinformatics (CBB) at Yale University.[21] The program was developed to address the challenges of interdisciplinary research related to the massive quantities of data produced by modern experimental technologies in the life sciences, aiming to bring together computer scientists, bioinformaticians, and experimental biomedical researchers. Creating such an environment for collaboration between these traditionally spatially and intellectually separated departments encourages new research ideas that incorporate the disparate perspectives of many different scientific disciplines. The PhD curriculum was designed to be equally cross-disciplinary and flexible enough to accommodate a population of trainees with heterogeneous backgrounds. Ideally, this type of program will produce students with a broader mindset and range of competencies more in line with emerging careers in computational science. Programs such as this can provide important lessons for structuring an applied component of a mathematics PhD program.

**Proposed changes.** We propose a similarly adaptable program structure for mathematics PhD programs that allows the unique interests and career aspirations of each student to influence the design of their graduate experience. These individualized needs can guide the relative emphasis of the major program components so that students are able to spend the most time on activities that will best encourage their professional development. As a transitional step to completely individualized training, Ph.D. programs could support several distinct "tracks", each with a structure tailored toward a particular career trajectory (see Table 2). For instance, since current programs are already strong in training students for traditional careers in academia, minimal modifications would be required to accommodate this career track. We also suggest a pedagogical track with a similar course load but a decreased research requirement for students hoping to pursue collegiate lecturer positions without substantial research components, allowing more time for a broader range of more in-depth teaching experiences.

A model similar to Yale's CBB program would be well-suited to students aiming for employment in applied math or industry -- the primary emphasis would be on gaining hands-on experience with the type of interdisciplinary research integral to such careers. A portion of the courses in this track should be external to the mathematics department (depicted as the lighter green coursework color in Table 2) to ensure that students obtain an adequate background in the applications of mathematical methods as well as the underlying theoretical frameworks. Ideally, this could also improve communication and integration with other disciplines outside the mathematics department. It should also be possible, or perhaps even required, for these students to integrate semester- or year-long internships in industry or government. This could be implemented through a network of professional opportunities curated by the program (as is currently done to some extent by the math program at New York University[25])as well as via existing infrastructure such as the BIG Math Network[26], Internship Network in the Mathematical Sciences (Inmas)[27], and AAU PhD Initiative[28]. Such opportunities provide valuable work experience, exposure to non-academic fields, and help students begin to build a professional network of more academically diverse peers than would typically be possible inside a mathematics department. Encouraging more students to pursue internships would benefit the institution in many ways as well. In addition to providing an additional mechanism of funding graduate students, it would expand and diversify the competencies of students beyond the areas of expertise of mathematics faculty, which could open the door to novel research ideas inspired by experiences in industry. Further, students would develop practical skills needed to do research earlier in their training, allowing them to have more time for independent work, successful grant applications, and publications, which would in turn benefit the research standing of the university in addition to the students.

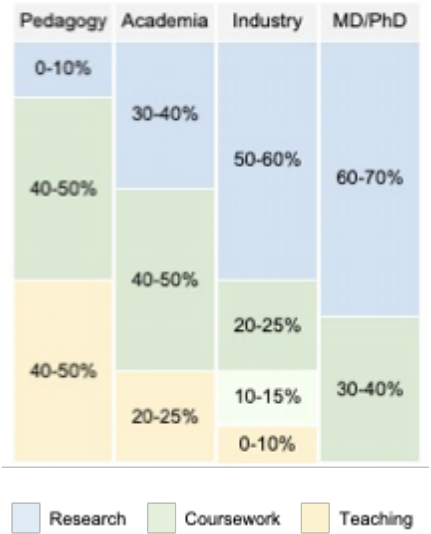

**Table 2.** Proposed composition of tailored graduate programs.

| Pedagogy | Academia | Industry | MD/PhD |
|---|---|---|---|
| 0-10% | 30-40% | 50-60% | 60-70% |
| 40-50% | 40-50% | 20-25% | |
| 40-50% | 20-25% | 10-15% | 30-40% |
| | | 0-10% | |

Legend: Research · Coursework · Teaching

Though it is rare for MD/PhD students to pursue their graduate training in mathematics, more might choose to do so if programs were better able to adapt to the unique requirements of a four-year PhD. We propose eliminating teaching requirements in this case since they are not needed to secure stipend funding. Given the accelerated timeline, we suggest course and departmental exam requirements could be lightened somewhat outside the chosen area of research focus. Courses completed through the medical school curriculum could serve a similar role to the external courses in the applied mathematics track to compensate. This decrease in teaching and course requirements would free additional time for MD/PhD students to carry out a level of research closer to that of their peers in less course-heavy graduate programs.

Finally, for tailored programs to be successful, high-quality professional guidance must be available and encouraged throughout the graduate experience. In the early years of the program, this should focus on educating students about the full breadth of careers open to them beyond positions in academia. Many incoming PhD students may lack the knowledge and experiences associated with different mathematical careers to make an informed choice about which track to choose upon entering the program. As such, faculty mentors will need to help these students explore their interests and choose program elements that will develop those interests into career prospects. Additionally, passage from one track to another should be flexible-- students may enter the program convinced they want to teach, but have a positive experience with an internship that piques their interest in an industry career, for example. As students progress, guidance can focus more on research-related challenges, and eventually should prepare students for the transition from graduate school into the next phase in their career.

**Potential roadblocks.** Of course, no program structure solves all problems, and the path to practical implementations of this type of student-centered PhD program would not be easy, requiring a host of nontrivial structural changes to existing programs. Novel challenges introduced by this program structure would need to be addressed, such as finding alternative means of financially supporting students whose teaching responsibilities have been reduced, compensating for the teaching burden displaced by such students, and ensuring departments are composed of faculty with sufficiently diverse backgrounds to provide high-quality research and career guidance to students in non-academic tracks. It is also important to consider the impact these changes could have on diversity, equity, and inclusion in math departments. These issues are all explored further in a companion article, "Implications of Student-Centric Graduate Programs in Mathematics", which we hope can initiate more open-ended discussions of these topics within the mathematical sciences community.

**References**

Alan Leshner and Layne Sherer (eds.), Graduate STEM Education for the 21st Century, National Academies Report, National Academies Press, Washington D.C., 2018.

12. XX

MD/PhD Programs

Case Studies

Mentoring

Rankings

Internships